\documentclass[12pt]{article}

\usepackage[authoryear]{natbib}
\usepackage{epsf,psfrag,amsmath,amsfonts,amssymb,bbm,epsfig}

\begin{document}

\title{Convergence of the empirical process in Mallows distance, with an application to bootstrap performance}
\author{Richard Samworth\thanks{Statistical Laboratory, 
DPMMS/CMS, University of Cambridge, Wilberforce Road, Cambridge CB3 0WB, UK. 
Fax: +44 1223 337956} \thanks{Email: \tt{rjs57@cam.ac.uk}} 
\and Oliver Johnson$^*$\thanks{Email: \tt{otj1000@cam.ac.uk}} }
\date{\today}
\maketitle

Running Title: Empirical process in Mallows distance



\newtheorem{theorem}{Theorem}[section]
\newtheorem{lemma}[theorem]{Lemma}
\newtheorem{proposition}[theorem]{Proposition}
\newtheorem{corollary}[theorem]{Corollary}
\newtheorem{conjecture}[theorem]{Conjecture}
\newtheorem{definition}[theorem]{Definition}
\newtheorem{example}[theorem]{Example}
\newtheorem{remark}[theorem]{Remark}
\newtheorem{condition}{Condition}
\newtheorem{main}{Theorem}
\newtheorem{assumption}[theorem]{Assumption}
\setlength{\parskip}{\parsep}
\setlength{\parindent}{0pt}

\def \outlineby #1#2#3{\vbox{\hrule\hbox{\vrule\kern #1% 
\vbox{\kern #2 #3\kern #2}\kern #1\vrule}\hrule}}%
\def \endbox {\outlineby{4pt}{4pt}{}}%
\newenvironment{proof}
{\noindent{\bf Proof\ }}{{\hfill \endbox
}\par\vskip2\parsep}
\newenvironment{pfof}[2]{\removelastskip\vspace{6pt}\noindent
 {\it Proof  #1.}~\rm#2}{\par\vspace{6pt}}

\hfuzz20pt

\newcommand{\Section}[1]{\setcounter{equation}{0} \section{#1}}
\newcommand{\var}{{\rm{Var\;}}}
\newcommand{\cov}{{\rm{Cov\;}}}
\newcommand{\tends}{\rightarrow \infty}
\newcommand{\ep}{{\mathbb {E}}}
\newcommand{\pr}{{\mathbb {P}}}
\newcommand{\re}{{\mathbb {R}}}
\newcommand{\I}{\mathbbm{1}}
\newcommand{\vc}[1]{{\bf {#1}}}
\newcommand{\conpr}{\stackrel{p}{\rightarrow}}
\newcommand{\cond}{\stackrel{d}{\rightarrow}}
\newcommand{\condo}{\stackrel{d^{\circ}}{\rightarrow}} 
\newcommand{\q}{{\cal{Q}}}
\newcommand{\blah}[1]{}

\begin{abstract}
We study the rate of convergence of the Mallows distance between the empirical distribution of a sample and the underlying population.  The surprising feature of our results is that the convergence rate is slower in the discrete case than in the absolutely continuous setting.  We show how the hazard function plays a significant role in these calculations.  As an application, we recall that the quantity studied provides an upper bound on the distance between the bootstrap distribution of a sample mean and its true sampling distribution.  Moreover, the convenient properties of the Mallows metric yield a straightforward lower bound, and therefore a relatively precise description of the asymptotic performance of the bootstrap in this problem. 
\end{abstract}

\section{Introduction and main results}
\label{section:intro}

Different problems in Probability and Statistics naturally lead to different choices of probability metric.  One such choice is the Mallows distance, also known as the Wasserstein or Kantorovich distance.  This metric has found extensive applications to a wide variety of fields; see \citet{Rachev1984} for a review.
\begin{definition}
\label{def:mallows}
For  $r  \geq  1$,  let  $\mathcal{F}_r$ denote the set of distribution functions $F$ satisfying $\int_{-\infty}^{\infty} |x|^r  \, dF(x) < \infty$.  For $F,G \in \mathcal{F}_r$, the Mallows metric $d_r(F,G)$ is defined by
\[
d_r(F,G) = \inf_{\mathcal{T}_{X,Y}} \bigl\{\mathbb{E}|X - Y|^r\bigr\}^{1/r},
\]
where  $\mathcal{T}_{X,Y}$ is the  set of  all joint  distributions of pairs of random variables $X$ and $Y$ whose marginal distributions are $F$  and $G$ respectively.  We also write $d_r(X,Y)$ for $d_r(F,G)$, where this will cause no confusion.
\end{definition}
The empirical process is a fundamental quantity of interest.  Though often not explicitly recognised as such, the Mallows distance between the empirical distribution and the underlying population has arisen in the work of several authors, including \citet{CsorgoHorvath1990} and \citet{delBarrioetal1999,delBarrioetal2000}.

Suppose $X_1, \ldots X_n$ are independent random variables, each having distribution function $F$ with mean $\mu$ and finite variance $\sigma^2 > 0$, and let
$\hat{F}_n$ denote the empirical distribution function of the sample, given by
\[
\hat{F}_n(x) = \frac{1}{n} \sum_{i=1}^n \I_{\{X_i \leq x \}}.
\]

The main contributions of this paper are threefold:
\begin{enumerate}
\item{
In Theorem~\ref{thm:limdisc}, we show that in the case of a discrete underlying population of finite support, $n^{1/4}d_2(\hat{F}_n,F)$ converges to an explicit nondegenerate limiting distribution.  We contrast this with the $n^{1/2}$ normalisation required by the previously cited authors in the absolutely continuous case.}
\item{In Section~\ref{section:tail}, we study the tail conditions required by \citet{CsorgoHorvath1990} and \citet{delBarrioetal2000} for the convergence of $n^{1/2}d_2(\hat{F}_n,F)$ in the absolutely continuous case.  In particular, by considering the hazard function, we show that one of the conditions of \citet{delBarrioetal2000} is redundant, and the statement of Theorem~2.1 of \citet{CsorgoHorvath1990} may be simplified.}
\item{Section~\ref{section:boot} is devoted to an application of these results.  We recall the calculation of \citet{ShaoTu1995}, showing that $d_2(\hat{F}_n,F)$ provides an upper bound on the Mallows distance between the bootstrap distribution of the sample mean and its true sampling distribution.  We give a straightforward lower bound on this latter quantity, yielding conditions under which the upper and lower bounds are of the same order.}
\end{enumerate}

\section{Convergence rates and limiting distributions for $d_r(\hat{F}_n,F)$}
\label{section:limit}

First, recall the following two lemmas about $d_r$, which are proved in \citet{Major1978} and \citet{BickelFreedman1981} respectively.
\begin{lemma}
For $F,G \in \mathcal{F}_r$, the infimum in Definition~\ref{def:mallows} is attained by the following  construction:  let  $U \sim  U(0,1)$, and set  $X = F^{-1}(U)$, $Y =  G^{-1}(U)$, where, for example, $F^{-1}(p) = \inf\{x \in \mathbb{R}:F(x) \geq p\}$.  Thus
\[
d_r(F,G) = \biggl(\int_0^1 |F^{-1}(p) - G^{-1}(p)|^r \, dp \biggr)^{1/r}.
\]
\end{lemma}
\begin{lemma} \label{lem:equivmetr}
If $(F_n)  \in \mathcal{F}$ and $F \in  \mathcal{F}$, then $d_r(F_n,F)
\rightarrow 0$  as $n  \rightarrow \infty$ if  and only if,  for every
bounded, continuous  function $g: \mathbb{R}  \rightarrow \mathbb{R}$,
we have both
\begin{enumerate}
\item{
$ \begin{displaystyle}
\lim_{n \rightarrow \infty} \int_{-\infty}^{\infty} g(x) \, dF_n(x) = \int_{-\infty}^{\infty} g(x) \, dF(x); 
\end{displaystyle}$
}
\item{
$ \begin{displaystyle}
\lim_{n \rightarrow \infty} \int_{-\infty}^{\infty} |x|^r \, dF_n(x) = \int_{-\infty}^{\infty} |x|^r dF(x).
\end{displaystyle}$
}
\end{enumerate}
Thus, convergence in the Mallows metric $d_r$ is equivalent to convergence in distribution together with convergence of the $r$th absolute moments.\end{lemma}

It follows immediately by Lemma \ref{lem:equivmetr} and the strong law of large numbers that $d_2(\hat{F}_n,F) \rightarrow 0$ almost surely as $n \rightarrow \infty$.  It is important to remark that although we can calculate a rate of convergence of the two parts above (that is, convergence in distribution and of the $r$th absolute moment), this will not help us find a rate of convergence of $d_2(\hat{F}_n, F)$, and we must use other techniques.  

When $F$ has a density $f$, \citet{CsorgoHorvath1990} and \citet{delBarrioetal2000} give conditions under which $n^{1/2}d_2(\hat{F}_n,F)$ converges to a nondegenerate limiting distribution.  A simple version of such results is given in Theorem~\ref{thm:rate} below.  However, these results do not cover the case of a discrete underlying population, which is studied later, in Theorem~\ref{thm:limdisc}.

In order to prove Theorem~\ref{thm:rate}, we need to introduce some notation.  Let $D = D[0,1]$ denote the space of left-continuous, real-valued functions on $[0,1]$ possessing right limits at each point.  We may equip $D$ with the uniform norm 
\[
\|x - y\|_{\infty} = \sup_{p \in [0,1]} |x(p) - y(p)|.
\]
A small complication arises from  the fact that the  normed space $(D,\|\cdot\|_{\infty})$  is non-separable, and the $\sigma$-algebra,  $\mathcal{D}$, generated  by the open  balls is strictly smaller than the Borel $\sigma$-algebra, $\mathcal{D}_{\mathrm{Borel}}$,  generated  by  the open  sets. This creates measurability problems, as explained in \citet{Chibisov1965}, which lead us to work with the space $(D,\mathcal{D},\|\cdot\|_{\infty})$.  A consequence of using the ball $\sigma$-algebra is that we must make a slight modification to the notion  of  weak convergence, in line with \citet{Billingsley1999}, p.~67:
\begin{definition}
If $(Y_n)_{n \geq 0}$ is a sequence of random elements of $(D,\mathcal{D},\|\cdot\|_{\infty})$, we write $Y_n \condo Y_0$ as $n \rightarrow \infty$ if
\[
\mathbb{E}\bigl(f(Y_n)\bigr) \rightarrow \mathbb{E}\bigl(f(Y_0)\bigr)
\]
as $n \rightarrow \infty$, for all bounded, continuous functions $f:D \rightarrow \mathbb{R}$ which are $\mathcal{D}$-measurable.   
\end{definition}
Throughout, $B = \bigl(B(p)\bigr)_{0 \leq p \leq 1}$ denotes a Brownian bridge; that is, a zero mean Gaussian process with 
\[
\mathrm{Cov}\bigl(B(p),B(q)\bigr) = p(1-q)
\]
for $p \leq q$.  For $p \in (0,1)$, let $\xi_p = F^{-1}(p) = \inf\{x \in \mathbb{R}: F(x) \geq p\}$.
\begin{theorem} \label{thm:rate}
Suppose that $F$ has a density $f$ such that $f(\xi_p)$ is positive and continuous for $p \in [0,1]$, and that the limits $\xi_0 = \lim_{p \searrow 0} \xi_p$ and $\xi_1 = \lim_{p \nearrow 1} \xi_p$ exist in $\mathbb{R}$.  Then
\[
n^{1/2}d_2(\hat{F}_n,F) \stackrel{d}{\rightarrow} \biggl(\int_0^1 B(p)^2(F^{-1})'(p)^2 \, dp \biggr)^{1/2} 
\]
as $n \rightarrow \infty$.
\end{theorem}
\begin{proof}
Theorem~1 on pp.~640--641 of \citet{ShorackWellner1986}, together with Corollary~1 on p.~48 of the same book, give that 
\begin{equation} 
\label{eq:convbr}
f(\xi_p)\,n^{1/2}\bigl(\hat{F}_n^{-1}(p) - \xi_p\bigr) \condo B(p)
\end{equation}
on $(D,\mathcal{D},\|\cdot\|_{\infty})$, as $n \rightarrow \infty$. Now, with probability one, $B$ belongs to the space $(C[0,1],\|\cdot\|_{\infty})$ of continuous real-valued functions on $[0,1]$ equipped with the uniform norm, and moreover this space is separable.  We can therefore apply the version of the continuous mapping theorem for $\condo$ convergence (Billingsley, 1999, pp.~67--68)  to a composition map $h(p) = h_2\bigl(h_1(p)\bigr)$ from $(D,\mathcal{D},\|\cdot\|_{\infty})$  to $\mathbb{R}$.  The individual maps $h_1:(D,\mathcal{D},\|\cdot\|_{\infty}) \rightarrow (D,\mathcal{D},\|\cdot\|_{\infty})$ and $h_2:(D,\mathcal{D},\|\cdot\|_{\infty}) \rightarrow \mathbb{R}$ are defined by
\[
h_1(G)(p) = \frac{G(p)^2}{f(\xi_p)^2} \quad \text{and} \quad h_2(G) = \biggl(\int_0^1 G(p) \, dp\biggr)^{1/2}.
\]
Observe that the continuity of $h_1$ follows from the fact that $f(\xi_p)$ attains its (positive) infimum for some $p \in [0,1]$.  We conclude that
\[
n^{1/2}d_2(\hat{F}_n,F) = \biggl(\int_0^1 n\bigl(\hat{F}_n^{-1}(p) - \xi_p\bigr)^2 \, dp\biggr)^{1/2} \condo 
\biggl(\int_0^1 \frac{B(p)^2}{f(\xi_p)^2} \, dp \biggr)^{1/2} 
\]
as $n \rightarrow \infty$.  The result follows on noting that any bounded, continuous function from $\mathbb{R}$ to $\mathbb{R}$ is (Borel) measurable.
\end{proof}
Theorem~\ref{thm:rate} is not the strongest possible; in particular, it is not necessary for the underlying density to have bounded support. \citet{CsorgoHorvath1990} and \citet{delBarrioetal2000} obtain the same conclusion with conditions which amount to control of the behaviour of the tails of $F$. For example, Equations (3.16) and (3.14) of \citet{delBarrioetal2000} consist of the following two
conditions:
\begin{condition} As $n \tends$,
\label{cond:3.16}
\[
n \int_{0}^{1/n} (\hat{F}_n^{-1}(p) - F^{-1}(p) )^2 \conpr 0 \ \ \text{and} \ \ n\int_{(n-1)/n}^1 \bigl(\hat{F}_n^{-1}(p) - F^{-1}(p)\bigr)^2 \, dp \conpr 0.
\]
\end{condition}
\begin{condition}
\label{cond:3.14}
\[
\int_0^1 \frac{p(1-p)}{f(F^{-1}(p))^2} \, dp = \int_{-\infty}^\infty \frac{F(x)\bigl(1-F(x)\bigr)}{f(x)} \, dx < \infty. 
\]
\end{condition}
A version of Theorem~2.1 of \citet{CsorgoHorvath1990} gives the following:
\begin{theorem}
\label{thm:Csorgo}
Suppose that $F$ has a density $f$ such that $f(\xi_p)$ is positive and continuous for $p \in (0,1)$, and monotone for $p$ sufficiently close to zero and one.  If Conditions \ref{cond:3.16} and \ref{cond:3.14} hold, then
\[
n^{1/2}d_2(\hat{F}_n,F) \stackrel{d}{\rightarrow} \biggl(\int_0^1 B(p)^2(F^{-1})'(p)^2 \, dp \biggr)^{1/2} 
\]
as $n \rightarrow \infty$.
\end{theorem}
Proposition~1 of \citet{delBarrioetal1999} verifies Condition~\ref{cond:3.16} when $F$ is the normal distribution.  We generalise this proposition in Section~\ref{section:tail}, by showing that Condition~\ref{cond:3.16} holds when the hazard function diverges in the tails.  In particular then, Condition~\ref{cond:3.16} will be seen to be a consequence of Condition~\ref{cond:3.14}. 

Now, we turn to the case of discrete random variables.  Note that in Theorem~\ref{thm:limdisc} below, we can understand the limiting distribution as the integral between 0 and 1 of a sum of delta functions at the points where $F^{-1}(\cdot)$ jumps.  Thinking of the limit in this way shows the analogy with the limits in Theorem~\ref{thm:rate} and Theorem~\ref{thm:Csorgo}.
\begin{theorem} \label{thm:limdisc}
Let $F$ be the distribution function corresponding to the probability mass function given by $\pr( X_1 = x_j) = p_j$, where $x_1 < x_2 < \ldots < x_m$ and $\sum_{j=1}^m p_j = 1$, with each $p_j > 0$.  For $j=1,\ldots,m$, let $q_j = p_1 + \ldots + p_j$.  Then
\[
n^{1/4} d_2 (\hat{F}_n,F) \cond \biggl(\, \sum_{j=1}^{m-1} |B(q_j)| (x_{j+1} - x_j)^2\biggr)^{1/2}
\]
as $n \rightarrow \infty$.
\end{theorem} 
\begin{proof}
\begin{figure}
\begin{center}
\psfrag{p1}[t][t]{$p$}
\psfrag{f1}[t][t]{$F^{-1}(p)$}
\psfrag{0}[b][t]{$0$}
\psfrag{q1}[b][t]{$q_1$}
\psfrag{q2}[b][t]{$q_2$}
\psfrag{q3}[b][t]{$q_3$}
\psfrag{x1}[t][b]{$x_1$}
\psfrag{x2}[t][b]{$x_2$}
\psfrag{x3}[t][b]{$x_3$}
\psfrag{x4}[t][b]{$x_4$}
\psfrag{1}[b][t]{$1$}
\epsfig{file=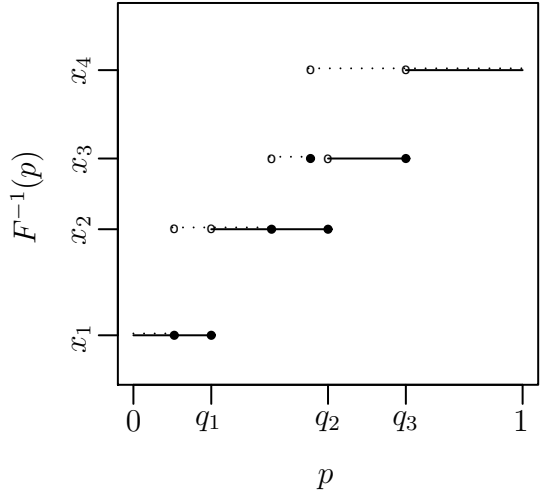}
\caption{\label{InverseDFGraph}Graphs of $F^{-1}(p)$ (solid) and $\hat{F}_n^{-1}(p)$ (dotted).  The dotted lines have been perturbed slightly to distinguish them from the solid ones.}
\end{center}
\end{figure}
See Figure~\ref{InverseDFGraph}.  For $j = 1,\ldots,m$, define $\hat{p}_j = n^{-1}\sum_{i=1}^n \I_{\{X_i = x_j\}}$ to be the empirical proportion of the sample taking the value $x_j$.  Let $\epsilon = \min_j p_j/3$, write $\hat{q}_j = \hat{p}_1 + \ldots + \hat{p}_j$, and define the event
\[
A = \{|\hat{q}_j - q_j| \leq \epsilon \ \text{for $j = 1,\ldots,m$}\}.
\]
Observe that, by the Dvoretsky--Keifer--Wolfowitz inequality \citep{Massart1990},
\[
\mathbb{P}(A^c) = \mathbb{P}\biggl(\bigcup_{j=1}^m |\hat{q}_j - q_j| > \epsilon\biggr) = \mathbb{P}\Bigl(\sup_{x \in \mathbb{R}}|\hat{F}_n(x) - F(x)| > \epsilon\Bigr) \leq 2e^{-2n\epsilon^2}.
\]
Thus, since $|\hat{F}_n^{-1}(p) - F^{-1}(p)| \leq x_m - x_1$,
\[
d_2^2(\hat{F}_n,F)\I_{A^c} \leq (x_m - x_1)^2\I_{A^c} = o_p(n^{-1/2})
\]
as $n \rightarrow \infty$.  On the other hand,
\[
d_2^2(\hat{F}_n,F)\I_A = \sum_{j=1}^{m-1} |\hat{q}_j - q_j| (x_{j+1} - x_j)^2.
\]
But $\bigl(n^{1/2}(\hat{q}_1 - q_1),\ldots,n^{1/2}(\hat{q}_m - q_m)\bigr) \cond N_m(0,\Sigma)$, where the asymptotic covariance matrix $\Sigma$ has entries $\Sigma_{ij} = q_i(1-q_j)$, for $1 \leq i \leq j \leq m$.  So by the continuous mapping theorem,
\[
n^{1/2}\sum_{j=1}^{m-1} |\hat{q}_j - q_j| (x_{j+1} - x_j)^2 \cond \sum_{j=1}^{m-1} |B(q_j)| (x_{j+1} - x_j)^2
\]
as $n \rightarrow \infty$.  Hence
\[
n^{1/2}d_2^2(\hat{F}_n,F) = n^{1/2}d_2^2(\hat{F}_n,F)\I_A + o_p(1) \cond \sum_{j=1}^{m-1} |B(q_j)| (x_{j+1} - x_j)^2
\]
as $n \rightarrow \infty$. 
\end{proof}
Theorems~\ref{thm:rate} and~\ref{thm:limdisc} extend to the case of general $d_r$.  Note the rates below coincide only when $r=1$, a case studied by \citet{delBarrioetal1999b}.  
\begin{corollary} \mbox{ }
\renewcommand{\theenumi}{\alph{enumi}}
\renewcommand{\labelenumi}{(\theenumi)}
\begin{enumerate}
\item{
Suppose that $F$ has a density $f$ such that $f(\xi_p)$ is positive and continuous for $p \in [0,1]$, and that the limits $\xi_0 = \lim_{p \searrow 0} \xi_p$ and $\xi_1 = \lim_{p \nearrow 1} \xi_p$ exist in~$\mathbb{R}$.  Then
\[
n^{1/2}d_r(\hat{F}_n,F) \cond \biggl(\int_0^1 |B(p)|^r (F^{-1})'(p)^r \, dp \biggr)^{1/2}
\]
as $n \rightarrow \infty$.}
\item{
Let $F$ be the distribution function corresponding to the probability mass function given by $\pr( X_1 = x_j) = p_j$, where $x_1 < x_2 < \ldots < x_m$ and $\sum_{j=1}^m p_j = 1$, with each $p_j > 0$.  Let $q_j = p_1 + \ldots + p_j$.  Then
\[
n^{1/(2r)} d_r (\hat{F}_n, F) \cond \biggl(\,\sum_{j=1}^{m-1} |B(q_j)|(x_{j+1}-x_j)^r\biggr)^{1/r}
\]
as $n \rightarrow \infty$.}
\end{enumerate}
\end{corollary} 
\begin{proof}
The proof of part $(a)$ follows that of Theorem~\ref{thm:rate}, the only difference being that the function $h_1$ should in this case be defined by
\[
h_1(G)(p) = \frac{|G(p)|^r}{f(\xi_p)^r}.
\]
Part $(b)$ mirrors the proof of Theorem~\ref{thm:limdisc}, noting that in general,
\[
n^{1/2}d_r^r(\hat{F}_n,F)\I_A = n^{1/2}\sum_{j=1}^{m-1}|\hat{q}_j - q_j|(x_{j+1} - x_j)^r,
\]
where the event $A$ was defined in the proof of Theorem~\ref{thm:limdisc}.
\end{proof}

\section{The hazard function and tail behaviour}
\label{section:tail}

In Section~\ref{section:limit}, we saw the importance of Conditions~1 and~2 concerning the tail behaviour of the distribution function $F$ in determining the rate of convergence to zero of $d_2(\hat{F}_n,F)$.  We now study the nature of these conditions in greater detail, using the notion of the hazard function. 
\begin{definition} 
For a random variable $X$ with distribution function $F$ and continuous density $f$, we define $\bar{F}(x) = \mathbb{P}(X > x)$, and the two-sided hazard function by
\[
h(x) = \left\{ \begin{array}{ll} f(x)/\bar{F}(x) & \mbox{if $x \geq \xi_{1/2}$,} \\
f(x)/F(x) & \mbox{if $x < \xi_{1/2}$.}
\end{array} \right.
\]
\end{definition}
The following theorem is the main result of this section, showing the redundancy of Condition \ref{cond:3.16}.  Its proof is deferred until after Lemma~\ref{lemma:mgf}.
\begin{theorem}
\label{thm:hazardtail}
If the hazard function satisfies $h(x) \rightarrow \infty$ as $|x| \rightarrow \infty$, then Condition \ref{cond:3.16} holds; that is,
\[
n\int_0^{1/n} \bigl(\hat{F}_n^{-1}(p) - F^{-1}(p)\bigr)^2 \, dp \conpr 0 \ \ \text{and} \ \ n\int_{(n-1)/n}^1 \bigl(\hat{F}_n^{-1}(p) - F^{-1}(p)\bigr)^2 \, dp \conpr 0
\]
as $n \rightarrow \infty$.  Moreover, Condition \ref{cond:3.14} implies 
Condition \ref{cond:3.16}. 
\end{theorem}
Let $X_t$ denote a random variable distributed as $X | X > t$, and consider $\var X_t$ (we set $\var X_t = 0$ if $\mathbb{P}(X > t) = 0$).  For example, if $X$ is an exponential random variable, then $X_t \sim X + t$, and we have $\var X_t = \var X$.  In general, this tail variance is closely related to the hazard function, as shown by the following lemma.
\begin{lemma} 
\label{lemma:haz}
Let $X$ be a random variable with hazard function $h(x)$.  If $t \geq \xi_{1/2}$, then
\[
\frac{1}{12 \sup_{x \geq t} h(x)^2} \leq \var X_t \leq \frac{4}{\inf_{x \geq t} h(x)^2},
\]
with a corresponding result for the left-hand tail.
\end{lemma}
\begin{proof}

\textbf{(a) Proof of the lower bound.}
Define the function
\[
k(t) = \frac{\sup_{x \geq t} f(x)}{\bar{F}(t)}.
\]
Note that the density of $X_t$ for $x \geq t$ is $f(x)/\bar{F}(t)$, which is at most $k(t)$.  By comparison with the $U[t,t+1/c]$ distribution, we see that if a random variable has density bounded above by $c$, then its variance is at least $1/(12 c^2)$, so $\var X_t \geq 1/\bigl(12 k(t)^2\bigr)$.  Moreover, if $x \geq t$, then $f(x)/\bar{F}(t) \leq f(x)/\bar{F}(x)$, from which we deduce that $k(t) \leq \sup_{x \geq t} h(x)$.

\textbf{(b) Proof of the upper bound.} In fact, we establish the stronger
conclusion that
\begin{equation} \label{eq:stronger}
\frac{ \ep \bigl( (X-t)^2 \I_{\{X>t\}} \bigr)}{\bar{F}(t)} \leq \frac{4}{\inf_{x \geq t} h(x)^2}. 
\end{equation}
We use arguments based on those which establish Poincar\'{e} inequalities, in \cite{BorovkovUtev1984} and \citet{Sysoeva1965}.  The simplest possible case of the main theorem of \citet{Sysoeva1965} gives a Hardy inequality, that if $G$ is a differentiable function with $G(t) = 0$, then
\begin{equation} \label{eq:sys}
\int_t^{\infty} f(x) G(x)^2 \, dx \leq 4 \int_t^\infty \frac{\bar{F}(x)^2}{f(x)} g(x)^2 \, dx = 4\int_t^\infty f(x) \frac{g(x)^2}{h(x)^2} \, dx,
\end{equation}
where $g(x) = G'(x)$, provided the integrals exist.

[For the sake of completeness, note that Equation (\ref{eq:sys}) can
be proved as follows. Integration by parts gives that
\begin{align*}
\int_t^\infty f(x) G(x)^2 \, dx & = 2 \int_t^\infty \bar{F}(x) G(x) g(x) \, dx \\
& \leq 2\biggl( \int_t^\infty f(x) G(x)^2 dx \biggr)^{1/2} \biggl(\int_t^\infty \frac{\bar{F}(x)^2}{f(x)} g(x)^2 dx \biggr)^{1/2},
\end{align*}
by Cauchy-Schwarz.]  Hence, from Equation (\ref{eq:sys}),
\[
\int_t^{\infty} f(x) G(x)^2 dx \leq \frac{4}{\inf_{x \geq t} h(x)^2} \int_t^\infty f(x) g(x)^2 dx,
\]
and choosing $G(x) = x - t$, Equation (\ref{eq:stronger}) follows.
\end{proof}
This shows that if $h(t)$ is bounded away from zero for large $|t|$, then $\var X_t < \infty$.  This is the case for the normal distribution, where the Mills ratio (Shorack and Wellner, 1986, p.~850) gives that $h(t) \geq |t|$ for all $t$.  On the other hand, if $h(t) \rightarrow 0$ as $|t| \rightarrow \infty$, then $\var X_t \rightarrow \infty$, as for the log-normal distribution, where $h(t) \sim t^{-1}\log t$ as $t \rightarrow \infty$.  Finally we remark that other tail variance behaviour is possible, in the case where $\liminf_{t \rightarrow \infty} h(t) < \limsup_{t \rightarrow \infty} h(t)$.  

As the example of the log-normal might suggest, results from reliability theory provide a link between the hazard function and the finiteness of the moment generating function in a neighbourhood of the origin.
\begin{lemma} 
\label{lemma:mgf}
Let $X$ be a random variable with hazard function $h(x)$.  If, for some $t \geq 0$, we have $\inf_{|x| \geq t} h(x) \geq c$, then $\mathbb{E}(e^{\theta X}) < \infty$ for $|\theta| < c$.
\end{lemma}
\begin{proof}
We prove the result for the right-hand tail.  Observe that for $x \geq t$,
\[
\log\biggl(\frac{\bar{F}(x)}{\bar{F}(t)}\biggr) = - \int_t^x h(y) \, dy,
\]
which can be seen by differentiating both sides with respect to $x$.  Thus  \begin{equation}
\label{eq:hazard}
\bar{F}(x) = \bar{F}(t) \exp \biggl(-\int_t^x h(y) \, dy \biggr).
\end{equation}
Now the right-hand side of Equation~(\ref{eq:hazard}) is bounded above by $e^{-cx}$, which is enough to guarantee that the moment generating function is finite for $\theta \in [0,c)$.
\end{proof}
Lemma~\ref{lemma:mgf} shows that if a random variable $X$ has a moment generating function which is infinite other than at the origin, then there exists a sequence $(x_n)$ tending to infinity such that $\lim_{n \tends} h(x_n) = 0$.  If, further, $X$ has a decreasing hazard function then we deduce that $\lim_{x \tends} h(x) =0$, and hence by Lemma~\ref{lemma:haz}, that $\lim_{t \tends} \var X_t= \infty$. 

We are now in a position to prove Theorem~\ref{thm:hazardtail}.

\begin{proof}{\textbf{of Theorem~\ref{thm:hazardtail}}}
We have
\begin{align}
n\int_{(n-1)/n}^1 \bigl(\hat{F}_n^{-1}&(p) - F^{-1}(p)\bigr)^2 \, dp \nonumber \\
&= n\int_{(n-1)/n}^1 \bigl(\hat{F}_n^{-1}(p) - \xi_{(n-1)/n} + \xi_{(n-1)/n} - F^{-1}(p)\bigr)^2 \, dp \nonumber \\ 
&\leq 2 \biggl( (X_{(n)} - \xi_{(n-1)/n})^2 + n \int_{\xi_{(n-1)/n}}^{\infty} (x - \xi_{(n-1)/n})^2 \, dF(x) \biggr).
\end{align}
The second term tends to zero by Equation~(\ref{eq:stronger}).  If $\xi_1 < \infty$, then the first term clearly converges in probability to zero.  On the other hand, if $\xi_1 = \infty$, then the assumed hypothesis on the hazard function combined with Equation~(\ref{eq:hazard}) imply that for any $t > 0$,
\[
\lim_{x \tends} \frac{ \bar{F}(t+x)}{\bar{F}(x)} = 0.
\]
These are is precisely the conditions under which \citet{Galambos1978}, Theorem 4.1.2, proves convergence in probability of $X_{(n)} - \xi_{(n-1)/n}$ to zero.  This completes the first part of the proof.

Note that Condition \ref{cond:3.14} can be restated in terms of the hazard function; that is, by considering the regions $x \in [\xi_{1/2},\infty)$ and $x \in (-\infty,\xi_{1/2})$ separately, we see that Condition \ref{cond:3.14} holds if and only if 
\[
\int_{-\infty}^{\infty} \frac{1}{h(x)} \, dx < \infty.
\]
Finiteness of this integral implies that $h(x) \rightarrow \infty$ as $|x| \tends$, and so by the argument above, Condition \ref{cond:3.16} holds.
\end{proof}
Note for the normal distribution that $h(x) \rightarrow \infty$, but, again by the Mills ratio, 
\[
\int_1^{\infty} \frac{1}{h(x)} \, dx \geq \int_1^{\infty} \Bigl(\frac{1}{x} - \frac{1}{x^3}\Bigr) \, dx = \infty.
\]
Thus Condition~\ref{cond:3.16} holds (implying Proposition~1 of \citet{delBarrioetal1999}), but not Condition~\ref{cond:3.14}.

Finally in this section, we give a partial converse to Theorem~\ref{thm:Csorgo}.
\begin{proposition}
If $\var X_t \rightarrow \infty$ as $t \rightarrow \infty$, then $n^{1/2}d_2(\hat{F}_n,F) \stackrel{a.s.}{\rightarrow} \infty$ as $n \rightarrow \infty$.
\end{proposition}
\begin{proof}
For $i = 1,\ldots,n$, let $a_i = n \int_{(i-1)/n}^{i/n} F^{-1}(p) \, dp$.  We decompose $nd_2^2(\hat{F}_n,F)$ into a random and deterministic part as follows:
\begin{align*}
n d_2^2(\hat{F}_n, F) & = n \sum_{i=1}^n \int_{(i-1)/n}^{i/n} \bigl( (X_{(i)} - a_i)^2 + (F^{-1}(p) - a_i)^2 \bigr) \, dp \\
&= \sum_{i=1}^n (X_{(i)} - a_i)^2 + n \sum_{i=1}^n \int_{(i-1)/n}^{i/n} (F^{-1}(p) - a_i)^2  \, dp \\
&\geq (X_{(n)} - a_n)^2 + n \int_{(n-1)/n}^1 (F^{-1}(p) - a_n)^2 \, dp \\
&= (X_{(n)} - a_n)^2 + n \int_{\xi_{(n-1)/n}}^{\infty}  (x - a_n)^2 \, dF(x) \\  
&= (X_{(n)} - a_n)^2 + \var X_{\xi_{(n-1)/n}} \, . 
\end{align*}
Ignoring the random part, the result is immediate.
\end{proof}

\section{Application to the bootstrap}
\label{section:boot}

The bootstrap was introduced into Statistics in the landmark paper of \citet{Efron1979}, and gives a very general technique for approximating the distributions of roots (i.e. functions of the sample and parameters of the underlying population) of interest.  The key idea for its use in practice is that of resampling; given a sample $X_1,\ldots,X_n$, we draw a further sample $X_1^*,\ldots,X_n^*$ uniformly at random with replacement from the original sample, and perform our calculations based on this resample.  The opportunity to repeat this resampling procedure enables the practitioner to mimic drawing additional samples from the original population.    

Historically, one of the great early triumphs in the analysis the bootstrap was the paper of \citet{Singh1981}.  One of his main results was to show that the bootstrap distribution of a normalised sample mean converges to its true sampling distribution at rate $O(n^{-1})$, provided the underlying population is non-lattice, and sufficiently many moments exist.  This improves on the $O(n^{-1/2})$ convergence rate of its natural competitor, namely normal approximation.  Singh's theorem was stated in terms of the supremum distance between the respective distribution functions.  The Mallows metric was first considered in the context of the bootstrap by \citet{BickelFreedman1981}.

Recall that $X_1,\ldots,X_n$ are independent random variables, each with distribution function $F$ with mean $\mu$ and finite variance $\sigma^2 > 0$, that $\hat{F}_n$ is the empirical distribution function of the sample and that $\bar{X}_n = n^{-1}\sum_{i=1}^n X_i$ is the sample mean.  A standard procedure for constructing confidence intervals for $\mu$ is to invert a probability statement concerning a root such as $n^{1/2}(\bar{X}_n - \mu)$, whose sampling distribution under $F$ we denote by $H_n(F)$.  Conditional on $X_1,\ldots,X_n$, let $X_1^*,\ldots,X_n^*$ be a resample; that is, an independent and identically distributed sample drawn from $\hat{F}_n$.  

The nonparametric bootstrap estimates the sampling  distribution  of  $n^{1/2}(\bar{X}_n - \mu)$  by that of $n^{1/2}(\bar{X}_n^* - \bar{X}_n)$, where $\bar{X}_n^* = n^{-1}\sum_{i=1}^n   X_i^*$.    In   other   words,   conditional   on $X_1,\ldots,X_n$,  we approximate  $H_n(F)$ by  $H_n(\hat{F}_n)$.  The properties of the Mallows distance $d_2$ make it suitable for studying the performance of the bootstrap approximation in this context.  The calculation  below follows \citet{ShaoTu1995}, and uses results proved in \citet{BickelFreedman1981}.
\begin{align}
d_2\bigl(H_n(\hat{F}_n),H_n(F)\bigr) &= d_2\biggl(\frac{1}{n^{1/2}}\sum_{i=1}^n (X_i^* - \bar{X}_n) \, , \, \frac{1}{n^{1/2}}\sum_{i=1}^n (X_i - \mu)\biggr) \nonumber \\
&\leq \frac{1}{n^{1/2}}\biggl(\sum_{i=1}^n d_2(X_i^* - \bar{X}_n \, , \, X_i 
- \mu)^2\biggr)^{1/2} \nonumber \\
&= d_2(X_1^* - \bar{X}_n \, , \, X_1 - \mu) \nonumber \\
&= \left( d_2(X_1^*, X_1)^2 - (\bar{X}_n - \mu)^2 \right)^{1/2} \nonumber \\
&\leq d_2(X_1^*,X_1) \nonumber \\
&= d_2(\hat{F}_n,F). \label{eq:Shao}
\end{align}
Thus, in particular, the distance between the distribution of the root of interest, $H_n(F)$, and its bootstrap approximation, $H_n(\hat{F}_n)$, is stochastically dominated by the distance between the true and empirical distributions.  The first inequality and the conditions under which we obtain equality are studied in \citet{JohnsonSamworth2004}.

Having studied the upper bound $d_2( \hat{F}_n, F)$ in previous sections, our final result gives a lower bound on the rate of convergence of $d_2\bigl(H_n(\hat{F}_n),H_n(F)\bigr)$.  Observe in particular that if the hypotheses of Theorem~\ref{thm:rate} and Proposition~\ref{prop:delta} are both satisfied, then the rates of convergence of the two bounds are identical. 
\begin{proposition} 
\label{prop:delta} \mbox{ } 

\renewcommand{\theenumi}{\alph{enumi}}
\renewcommand{\labelenumi}{(\theenumi)}
\begin{enumerate}
\item If \, $\mathbb{E}X_1^4 < \infty$, then $n^{1/2}d_2\bigl(H_n(\hat{F}_n),H_n(F)\bigr)$ is bounded away from zero in probability; in fact, for every $\epsilon > 0$,
\[
\mathbb{P}\bigl\{n^{1/2}d_2\bigl(H_n(\hat{F}_n),H_n(F)\bigr) > \epsilon \bigr\} \nrightarrow 0
\]
as $n \rightarrow \infty$.
\item Fix $\delta \in (0,2)$.  If $\mathbb{E}|X_1|^{2+\delta} = \infty$, then, with probability one, \linebreak $n^{\delta/(2+\delta)}d_2\bigl(H_n(\hat{F}_n),H_n(F)\bigr)$ does not tend to zero; that is
\[
\mathbb{P}\bigl\{n^{\delta/(2+\delta)}d_2\bigl(H_n(\hat{F}_n),H_n(F)\bigr) \rightarrow 0\bigr\} = 0.
\]
\end{enumerate}
\end{proposition}
\begin{proof}
For any $r \geq 1$ and $F_X,F_Y \in \mathcal{F}_r$, there exist random variables $(X,Y)$ having marginal distribution functions $F_X$ and $F_Y$ respectively, and such that
\begin{equation} 
d_r(F_X,F_Y) = \bigl(\ep |X-Y|^r \bigr)^{1/r} \geq \bigl| (\ep |X|^r)^{1/r}  - (\ep |Y|^r)^{1/r}\bigr|.
\end{equation}
This means that 
\begin{align*}
d_2\bigl(H_n(\hat{F}_n),H_n(F)\bigr) &= d_2\biggl(\frac{1}{n^{1/2}}\sum_{i=1}^n (X_i^* - \bar{X}_n) \, , \, \frac{1}{n^{1/2}}\sum_{i=1}^n (X_i - \mu)\biggr) \\ 
&\geq |s - \sigma|
\end{align*}
where $s^2 = n^{-1}\sum (X_i - \bar{X}_n)^2$ is the sample variance.  Without loss of generality, suppose that $\mu = \ep X_1 = 0$.

To prove $(a)$, observe that if $\ep (X_1^4) < \infty$, then $n^{1/2}|s - \sigma|$ is bounded away from zero in probability; in fact, $n^{1/2}(s - \sigma) \cond N\bigl(0,(\ep X_1^4 - \sigma^4)/(4 \sigma^2)\bigr)$ (Serfling, 1980, p.~119). 

To prove $(b)$, note that 
\begin{align*}
n^{\delta/(2 + \delta)}(s^2 - \sigma^2) &= \frac{1}{n^{2/(2+\delta)}}\sum_{i=1}^n \{(X_i - \bar{X}_n)^2 - \sigma^2\} \\
&= \frac{1}{n^{2/(2+\delta)}}\sum_{i=1}^n (X_i^2 - \sigma^2) - n^{\delta/(2+\delta)}\bar{X}_n^2.
\end{align*}
Now, applying the law of the iterated logarithm gives that 
\[
n^{\delta/(2+\delta)}\bar{X}_n^2 = O(n^{-2/(2+\delta)}\log \log n) \quad \text{a.s..}
\]
Moreover, since $\ep|X_1|^{2+\delta} = \infty$, by the converse to the Marcinkiewicz-Zygmund strong law of large numbers (cf. \citet{Loeve1977}, p.~255),
\[
\mathbb{P}\biggl(\frac{1}{n^{2/(2+\delta)}}\sum_{i=1}^n (X_i^2 - \sigma^2) \rightarrow 0\biggr) < 1.
\]
Hence, by Kolmogorov's zero-one law, and since $s^2 - \sigma^2 = (s-\sigma)(s + \sigma)$,
\[
\mathbb{P}\bigl(n^{\delta/(2+\delta)}|s - \sigma| \rightarrow 0\bigr) = 0.
\]
\end{proof}
Note that the strength of the conclusion of part $(b)$ of Proposition~\ref{prop:delta} increases as $\delta$ decreases.  

\def\polhk#1{\setbox0=\hbox{#1}{\ooalign{\hidewidth
  \lower1.5ex\hbox{`}\hidewidth\crcr\unhbox0}}}

\end{document}